\documentclass{amsart} 
\usepackage{latexsym} 
\usepackage{amssymb} 
\usepackage{amsmath} 

\begin{document}

%%%%%%%%%%%%%%%%%%%%%%Definitions%%%%%%%%%%%%%%%%%%%%%%%%%%%%%%%%%%%%%%%%%%% 

\newtheorem{theorem}{Theorem} 
\newtheorem{problem}{Problem} 
\newtheorem{definition}{Definition} 
\newtheorem{lemma}{Lemma} 
\newtheorem{proposition}{Proposition} 
\newtheorem{corollary}{Corollary} 
\newtheorem{example}{Example} 
\newtheorem{conjecture}{Conjecture} 
\newtheorem{algorithm}{Algorithm} 
\newtheorem{exercise}{Exercise} 
\newtheorem{remarkk}{Remark} 
 
\newcommand{\be}{\begin{equation}} 
\newcommand{\ee}{\end{equation}} 
\newcommand{\bea}{\begin{eqnarray}} 
\newcommand{\eea}{\end{eqnarray}} 
\newcommand{\beq}[1]{\begin{equation}\label{#1}} 
\newcommand{\eeq}{\end{equation}} 
\newcommand{\beqn}[1]{\begin{eqnarray}\label{#1}} 
\newcommand{\eeqn}{\end{eqnarray}} 
\newcommand{\beaa}{\begin{eqnarray*}} 
\newcommand{\eeaa}{\end{eqnarray*}} 
\newcommand{\req}[1]{(\ref{#1})} 
 
\newcommand{\lip}{\langle} 
\newcommand{\rip}{\rangle} 

\newcommand{\uu}{\underline} 
\newcommand{\oo}{\overline} 
\newcommand{\La}{\Lambda} 
\newcommand{\la}{\lambda} 
\newcommand{\eps}{\varepsilon} 
\newcommand{\om}{\omega} 
\newcommand{\Om}{\Omega} 
\newcommand{\ga}{\gamma} 
\newcommand{\rrr}{{\Bigr)}} 
\newcommand{\qqq}{{\Bigl\|}} 
 
\newcommand{\dint}{\displaystyle\int} 
\newcommand{\dsum}{\displaystyle\sum} 
\newcommand{\dfr}{\displaystyle\frac} 
\newcommand{\bige}{\mbox{\Large\it e}} 
\newcommand{\integers}{{\Bbb Z}} 
\newcommand{\rationals}{{\Bbb Q}} 
\newcommand{\reals}{{\rm I\!R}} 
\newcommand{\realsd}{\reals^d} 
\newcommand{\realsn}{\reals^n} 
\newcommand{\NN}{{\rm I\!N}} 
\newcommand{\DD}{{\rm I\!D}} 
\newcommand{\degree}{{\scriptscriptstyle \circ }} 
\newcommand{\dfn}{\stackrel{\triangle}{=}} 
\def\complex{\mathop{\raise .45ex\hbox{${\bf\scriptstyle{|}}$} 
     \kern -0.40em {\rm \textstyle{C}}}\nolimits} 
\def\hilbert{\mathop{\raise .21ex\hbox{$\bigcirc$}}\kern -1.005em {\rm\textstyle{H}}} %Hilbert space 
\newcommand{\RAISE}{{\:\raisebox{.6ex}{$\scriptstyle{>}$}\raisebox{-.3ex} 
           {$\scriptstyle{\!\!\!\!\!<}\:$}}} % >< one above each other 
 
\newcommand{\hh}{{\:\raisebox{1.8ex}{$\scriptstyle{\degree}$}\raisebox{.0ex} 
           {$\textstyle{\!\!\!\! H}$}}} 

\newcommand{\OO}{\won} 
\newcommand{\calA}{{\mathcal A}} 
\newcommand{\calB}{{\cal B}} 
\newcommand{\calC}{{\cal C}} 
\newcommand{\calD}{{\cal D}} 
\newcommand{\calE}{{\cal E}} 
\newcommand{\calF}{{\mathcal F}} 
\newcommand{\calG}{{\cal G}} 
\newcommand{\calH}{{\cal H}} 
\newcommand{\calK}{{\cal K}} 
\newcommand{\calL}{{\mathcal L}} 
\newcommand{\calM}{{\mathcal M}} 
\newcommand{\calO}{{\cal O}} 
\newcommand{\calP}{{\cal P}} 
\newcommand{\calU}{{\mathcal U}} 
\newcommand{\calX}{{\cal X}} 
\newcommand{\calXX}{{\cal X\mbox{\raisebox{.3ex}{$\!\!\!\!\!-$}}}} 
\newcommand{\calXXX}{{\cal X\!\!\!\!\!-}} 
\newcommand{\gi}{{\raisebox{.0ex}{$\scriptscriptstyle{\cal X}$} 
\raisebox{.1ex} {$\scriptstyle{\!\!\!\!-}\:$}}} 
\newcommand{\intsim}{\int_0^1\!\!\!\!\!\!\!\!\!\sim} 
\newcommand{\intsimt}{\int_0^t\!\!\!\!\!\!\!\!\!\sim} 
\newcommand{\pp}{{\partial}} 
\newcommand{\al}{{\alpha}} 
\newcommand{\sB}{{\cal B}} 
\newcommand{\sL}{{\cal L}} 
\newcommand{\sF}{{\cal F}} 
\newcommand{\sE}{{\cal E}} 
\newcommand{\sX}{{\cal X}} 
\newcommand{\R}{{\rm I\!R}} 
\renewcommand{\L}{{\rm I\!L}} 
\newcommand{\vp}{\varphi} 
\newcommand{\N}{{\rm I\!N}} 
\def\ooo{\lip} 
\def\ccc{\rip} 
\newcommand{\ot}{\hat\otimes} 
\newcommand{\rP}{{\Bbb P}} 
\newcommand{\bfcdot}{{\mbox{\boldmath$\cdot$}}} 
 
\renewcommand{\varrho}{{\ell}} 
\newcommand{\dett}{{\textstyle{\det_2}}} 
\newcommand{\sign}{{\mbox{\rm sign}}} 
\newcommand{\TE}{{\rm TE}} 
\newcommand{\TA}{{\rm TA}} 
\newcommand{\E}{{\rm E\,}} 
\newcommand{\won}{{\mbox{\bf 1}}} 
\newcommand{\Lebn}{{\rm Leb}_n} 
\newcommand{\Prob}{{\rm Prob\,}} 
\newcommand{\sinc}{{\rm sinc\,}} 
\newcommand{\ctg}{{\rm ctg\,}} 
\newcommand{\loc}{{\rm loc}} 
\newcommand{\trace}{{\,\,\rm trace\,\,}} 
\newcommand{\Dom}{{\rm Dom}} 
\newcommand{\ifff}{\mbox{\ if and only if\ }} 
\newcommand{\nproof}{\noindent {\bf Proof:\ }} 
\newcommand{\remark}{\noindent {\bf Remark:\ }} 
\newcommand{\remarks}{\noindent {\bf Remarks:\ }} 
\newcommand{\note}{\noindent {\bf Note:\ }}

\newcommand{\boldx}{{\bf x}} 
\newcommand{\boldX}{{\bf X}} 
\newcommand{\boldy}{{\bf y}} 
\newcommand{\boldR}{{\bf R}} 
\newcommand{\uux}{\uu{x}} 
\newcommand{\uuY}{\uu{Y}} 
 
\newcommand{\limn}{\lim_{n \rightarrow \infty}} 
\newcommand{\limN}{\lim_{N \rightarrow \infty}} 
\newcommand{\limr}{\lim_{r \rightarrow \infty}} 
\newcommand{\limd}{\lim_{\delta \rightarrow \infty}} 
\newcommand{\limM}{\lim_{M \rightarrow \infty}} 
\newcommand{\limsupn}{\limsup_{n \rightarrow \infty}} 
 
\newcommand{\ra}{ \rightarrow }

\newcommand{\ARROW}[1] 
  {\begin{array}[t]{c}  \longrightarrow \\[-0.2cm] \textstyle{#1} \end{array} } 
 
\newcommand{\AR} 
 {\begin{array}[t]{c} 
  \longrightarrow \\[-0.3cm] 
  \scriptstyle {n\rightarrow \infty} 
  \end{array}} 
 
\newcommand{\pile}[2] 
  {\left( \begin{array}{c}  {#1}\\[-0.2cm] {#2} \end{array} \right) } 
 
\newcommand{\floor}[1]{\left\lfloor #1 \right\rfloor} 
 
%for doing boldface subscripts etc., e.g. $G_{\mmbox{\boldx}}$ 
\newcommand{\mmbox}[1]{\mbox{\scriptsize{#1}}} 
 
%fraction with round brackets 
\newcommand{\ffrac}[2] 
  {\left( \frac{#1}{#2} \right)} 
 
\newcommand{\one}{\frac{1}{n}\:} 
\newcommand{\half}{\frac{1}{2}\:} 
 
\def\le{\leq} 
\def\ge{\geq} 
\def\lt{<} 
\def\gt{>} 
 
%qed 
\def\squarebox#1{\hbox to #1{\hfill\vbox to #1{\vfill}}} 
\newcommand{\nqed}{\hspace*{\fill} 
           \vbox{\hrule\hbox{\vrule\squarebox{.667em}\vrule}\hrule}\bigskip} 
%%%%%%%%%%%%%%%%%%%%%%%%%%%%%%%%%%%%%%%%%%%%%%%%%%%%%%%%%%%%%%%%%%%%%%%%%%%% 
 
\title{Regularity of the conditional expectations with respect to
  signal to noise  ratio}

\author{ A. S. \"Ust\"unel\\Tot passa, per\`o segueix sent l'amistat} 
%\date{ } 
\maketitle 
\noindent 
{\bf Abstract:}{\small{ Let $(W,H,\mu)$ be the classical  Wiener space, assume 
that $U_\la=I_W+u_\la$ is an adapted perturbation of identity where the
perturbation $u_\la$ is an  $H$-valued map, defined up to
$\mu$-equivalence classes, such that its Lebesgue density $s\to
\dot{u}_\la(s)$ is almost surely  adapted
to the canonical  filtration of the Wiener space and depending
measurably on a real parameter $\la$. Assuming some   regularity  for
$u_\la$, its Sobolev derivative and integrability of the divergence of
the resolvent operator of its Sobolev derivative, we prove the almost sure and
$L^p$-regularity w.r. to $\la$  of the estimation  $E[\dot{u}_\la(s)|\calU_\la(s)]$
and more generally of the conditional expectations of the type
$E[F\mid\calU_\la(s)]$ for nice Wiener functionals, where 
 $(\calU_\la(s),s\in [0,1])$ is the the filtration  which is  generated by $U_\la$.
These results are applied to prove the invertibility of the adapted
perturbations of identity, hence to prove the strong  existence and
uniqueness of functional SDE's;  convexity of the entropy and the
quadratic estimation error and finally to the information theory.}}\\ 

\vspace{0.5cm} 

\tableofcontents
\section{\bf{Introduction} }

\noindent
The Malliavin calculus studies  the  regularity of the laws of
the random variables (functionals) defined on a Winer space (abstract
or classical) with values in finite dimensional Euclidean spaces (more
generally manifolds) using a variational calculus in the direction of
the underlying quasi-invariance space, called the Cameron-Martin
space. Although its  efficiency is  globally recognized by  now, for the maps taking values in the infinite
dimensional spaces the Malliavin calculus does not apply as easily as  in the
finite dimensional case due to the absence of the Lebesgue measure and
even the problem itself 
needs to be defined. For instance, there is a notion called signal to
noise  ratio which finds its roots in engineering which requires
regularity of infinite dimensional objects with respect to finite
dimensional parameters (cf.\cite{D,GSV,G-Y,KZZ,M-Z,P}). Let us explain
the problem along  its general lines  briefly: imagine  a communication
channel of the form $y=\sqrt{\la}x+w$, where $x$ denotes the emitted
signal and $w$ is a noise which corrupts the communications. The
problem of estimation of the signal $x$ from the data generated $y$ is
studied since the early beginnings of the electrical engineering. One
of the main problems dealt with is the behavior of the $L^2$-error of
the estimation w.r. to the signal to noise ratio $\la$. This requires
elementary probability when $x$ and $w$ are independent finite
dimensional variables, though it gives important results for
engineers. In particular, it has been recently realized that (cf.
\cite{GSV,Z}), in this linear model with $w$ being Gaussian,  the derivative of the mutual information between $x$ and
$y$ w.r. to $\la$ equals to the half of the mean quadratic error of estimation. 
 The infinite dimensional case is more tricky and requires
already the techniques of Wiener space analysis and the Malliavin
calculus (cf. \cite{Z}). The situation is much more complicated in the
case where the signal is correlated to
the noise; in fact we need the $\la$-regularity  of the conditional expectations w. r. to
the filtration generated by $y$, which is, at first sight,  clearly outside the scope of
the Malliavin calculus. 

In this paper we study the  generalization of the problem mentioned
above. Namely assume that we are given, in the setting of a classical Wiener space,
denoted as $(W,H,\mu)$, 
a signal which is of the form of an adapted perturbation of identity:
$$
U_\la(t,w)=W_t(w)+\int_0^t\dot{u}_\la(s,w)ds\,,
$$
where $(W_t,t\in[0,1])$ is the canonical Wiener process, $\dot{u}_\la$ is
an element of $L^2(ds\times d\mu)$ which is adapted to the Brownian
filtration $ds$-almost surely and $\la$ is a real parameter. Let
$\calU_\la(t)$ be the sigma algebra generated by $(U_\la(s),\,s\leq
t)$. What can we say about the regularity, i.e., continuity and/or
differentiability w.r. to $\la$,  of the functionals of the form 
$\la\to E[F\mid \calU_\la(t)]$ and $\la\to E[F\mid U_\la=w]$ 
(the latter denotes the disintegration) given various regularity
assumptions about the map $\la\to \dot{u}_\la$, like 
differentiability of it or its $H$-Sobolev derivatives w.r. to $\la$?
We prove that the answer to these questions depend essentially on the
behavior of the random  resolvent operator $(I_H+\nabla u_\la)^{-1}$,
where $\nabla u_\la$ denotes the Sobolev derivative of $u_\la$, which
is a quasi-nilpotent Hilbert-Schmidt operator, hence its resolvent
exists always.  More precisely  we prove that if the functional
\begin{equation}
\label{cond-1}
(1+\rho(-\delta u_\la)\delta\left((I_H+\nabla u_\la)^{-1}\frac{d}{d\la}u_\la\right)
\end{equation}
is in $L^1(d\la\times d\mu,[0,M]\times W)$ for some $M>0$, where
$\delta$ denotes the Gaussian divergence and $\rho(-\delta u)$ is the
Girsanov-Wick exponential corresponding to the stochastic integral
$\delta u$ (cf. the next section),  then the
map $\la\to L_\la$ is absolutely continuous almost surely where
$L_\la$ is the Radon-Nikodym derivative of $U_\la\mu$ w.r. to $\mu$
and we can calculate its derivative explicitly. This observation
follows from some variational calculus and from  the Malliavin
calculus. The iteration of the hypothesis (\ref{cond-1}) by replacing
$\delta((I_H+\nabla u_\la)^{-1}\frac{d}{d\la}u_\la)$
with its $\la$-derivatives permits us to prove the higher order  differentiability
of the above conditional expectations w.r. to $\la$ and these results
are exposed in  Section \ref{S-Bas}. In Section \ref{S-Inv}, we
give applications of these results to show the almost sure
invertibility of the adapted perturbations of the identity, which is
equivalent to the strong  existence and uniqueness results of the
(functional) stochastic differential equations. In Section
\ref{S-Ent}, we apply the results of Section \ref{S-Bas} to calculate
the derivatives of the relative  entropy of $U_\la\mu$ w.r. to $\mu$
in the general case, i.e., we do not suppose the a.s. invertibility of
$U_\la$, which demands the calculation of the derivatives of the
non-trivial conditional expectations. Some results are also given for
the derivative of the quadratic error in the  case of anticipative
estimation as well as the relations to the  Monge-Kantorovich measure
transportation theory and the Monge-Amp\`ere equation. In Section
\ref{S-Info}, we generalize the celebrated result about the relation
between the mutual information and the mean quadratic error
(cf. \cite{D,G-Y,KZZ}) in the following way: we suppress the  hypothesis of
independence between the signal and the noise as well as the almost
sure invertibility of the observation for fixed  exterior 
parameter of the signal. With the help of the results of Section
\ref{S-Bas}, the calculations of the first and second
order  derivatives of the mutual information w.r. to the ratio
parameter $\la$ are also  given.

\section{\bf{Preliminaries and notation}}
\label{S-Pr}
\label{preliminaries}
\noindent
Let $W$ be the classical Wiener  space $C([0,T],\R^n)$  with  the Wiener
measure $\mu$. The corresponding Cameron-Martin space is denoted by $H$. Recall that the
injection $H\hookrightarrow W$ is compact and its adjoint is the
natural injection $W^\star\hookrightarrow H^\star\subset
L^2(\mu)$. Since the image of  $\mu$ under the mappings $w\to
w+h,\,h\in H$ is 
equivalent to $\mu$, the G\^ateaux  derivative in the  $H$ direction of the
random variables is a closable operator on $L^p(\mu)$-spaces and  this
closure is  denoted by $\nabla$ and called the Sobolev derivative (on
the Wiener space) cf., for example
\cite{ASU, ASU-1}. The corresponding Sobolev spaces consisting of 
(the equivalence classes) of  real-valued  random variables
will be denoted as $\DD_{p,k}$, where $k\in \NN$ is the order of
differentiability and $p>1$ is the order of integrability. If the
random variables are with values in some separable Hilbert space, say
$\Phi$, then we shall define similarly the corresponding Sobolev
spaces and they are denoted as $\DD_{p,k}(\Phi)$, $p>1,\,k\in
\NN$. Since $\nabla:\DD_{p,k}\to\DD_{p,k-1}(H)$ is a continuous and
linear operator its adjoint is a well-defined operator which we
represent by $\delta$.  A very important feature in the theory is that
$\delta$ coincides with the It\^o integral of the Lebesgue density of the adapted elements of
$\DD_{p,k}(H)$ (cf.\cite{ASU,ASU-1}).

For any $t\geq 0$ and measurable $f:W\to \reals_+$, we note by
$$
P_tf(x)=\int_Wf\left(e^{-t}x+\sqrt{1-e^{-2t}}y\right)\mu(dy)\,,
$$
it is well-known that $(P_t,t\in \reals_+)$ is a hypercontractive
semigroup on $L^p(\mu),p>1$,  which is called the Ornstein-Uhlenbeck
semigroup (cf.\cite{ASU,ASU-1}). Its infinitesimal generator is denoted
by $-\calL$ and we call $\calL$ the Ornstein-Uhlenbeck operator
(sometimes called the number operator by the physicists). The
norms defined by
\begin{equation}
\label{norm}
\|\phi\|_{p,k}=\|(I+\calL)^{k/2}\phi\|_{L^p(\mu)}
\end{equation}
are equivalent to the norms defined by the iterates of the  Sobolev
derivative $\nabla$. This observation permits us to identify the duals
of the space $\DD_{p,k}(\Phi);p>1,\,k\in\NN$ by $\DD_{q,-k}(\Phi')$,
with $q^{-1}=1-p^{-1}$,
where the latter  space is defined by replacing $k$ in (\ref{norm}) by
$-k$, this gives us the distribution spaces on the Wiener space $W$
(in fact we can take as $k$ any real number). An easy calculation
shows that, formally, $\delta\circ \nabla=\calL$, and this permits us
to extend the  divergence and the derivative  operators to the
distributions as linear,  continuous operators. In fact
$\delta:\DD_{q,k}(H\otimes \Phi)\to \DD_{q,k-1}(\Phi)$ and
$\nabla:\DD_{q,k}(\Phi)\to\DD_{q,k-1}(H\otimes \Phi)$ continuously, for
any $q>1$ and $k\in \reals$, where $H\otimes \Phi$ denotes the
completed Hilbert-Schmidt tensor product (cf., for instance
\cite{ASU,ASU-1,BOOK}). We shall denote by $\DD(\Phi)$ and $\DD'(\Phi)$
respectively the sets
$$
\DD(\Phi)=\bigcap_{p>1,k\in \N}\DD_{p,k}(\Phi)\,,
$$
and 
$$
\DD'(\Phi)=\bigcup_{p>1,k\in \N}\DD_{p,-k}(\Phi)\,,
$$
where the former is equipped with the projective and the latter is
equipped with the inductive limit topologies.

Let us denote by $(W_t,t\in[0,1])$ the coordinate map on $W$ which is
the canonical  Brownian motion (or Wiener process) under the Wiener
measure, let $(\calF_t,t\in[0,1])$ be its completed filtration. The
elements  of $L^2(\mu,H)=\DD_{2,0}(H)$ such that $w\to\dot{u}(s,w)$
are $ds$-a.s. $\calF_S$ measurable will be noted as $L^2_a(\mu,H)$ or
$\DD^a_{2,0}(H)$. $L^0_a(\mu,H)$ is defined similarly (under the
convergence in probability). Let $U:W\to W$ be defined as $U=I_W+u$
with some $u\in L^0_a(\mu,H)$, we say that $U$ is $\mu$-almost surely
invertible if there exists some $V:W\to W$ such that $V\mu\ll\mu$ and
that 
$$
\mu\left\{w:U\circ V(w)=V\circ U(w)=w\right\}=1\,.
$$
The following results are proved with various extensions  in
\cite{ASU-2,ASU-3,ASU-4}:
\begin{theorem}
\label{thm-0}
Assume that $u\in L^0_a(\mu,H)$, let $L$ be the Radon-Nikodym  density
of $U\mu=(I_W+u)\mu$ w.r. to $\mu$, where $U\mu$ denotes the image
(push forward) of $\mu$ under the map $U$. Then we have 
\begin{enumerate}
\item
$$
E[L\log L]\leq
\frac{1}{2}\|u\|_{L^2(\mu,H)}^2=\frac{1}{2}E\int_0^1|\dot{u}_s|^2ds\,.
$$
\item Assume that $E[\rho(-\delta u)]=1$, then we have the equality:
\begin{equation}
\label{eqlty}
E[L\log L]=\frac{1}{2}\|u\|_{L^2(\mu,H)}^2
\end{equation}
if and only if $U$ is almost surely invertible and its inverse can be
written as $V=I_W+v$, with $v\in L_a^0(\mu,H)$. 
\item
Assume that $E[L\log L-\log L]<\infty$ and the equality (\ref{eqlty})
holds, then $U$  is again  almost surely invertible and its inverse can be
written as $V=I_W+v$, with $v\in L_a^0(\mu,H)$.
\end{enumerate}
\end{theorem}

\noindent
The following result gives the relation between the entropy and the
estimation ( cf. \cite{ASU-2} for the proof):
\begin{theorem}
\label{thm-00}
Assume that $u\in L^2_a(\mu,H)$, let $L$ be the Radon-Nikodym  density
of $U\mu=(I_W+u)\mu$ w.r. to $\mu$, where $U\mu$ denotes the image
(push forward) of $\mu$ under the map $U$ and let $(\calU_t,t\in
[0,1])$ be the filtration generated by $(t,w)\to U(t,w)$. Assume that
$E[\rho(-\delta u)]=1$.  Then we have 
\begin{itemize}
\item
$$
E[L\log L]= \frac{1}{2}E\int_0^1|E[\dot{u}_s\mid \calU_s]|^2ds\,.
$$
\item 
$$
L\circ U\,E[\rho(-\delta u)|U]=1
$$
$\mu$-almost surely.
\end{itemize}
\end{theorem}

\section{\bf{Basic results}}
\label{S-Bas}
\noindent
Let $(W,H,\mu)$ be the classical Wiener space, i.e.,
$W=C_0([0,1],\R^d),\,H=H^1([0,1],\R^d)$ and $\mu$ is the Wiener
measure under which the evaluation map at $t\in [0,1]$ is a Brownian motion. Assume that $U_\la:W\to W$ is defined as
$$
U_\la(t,w)=W_t(w)+\int_0^t\dot{u}_\la(s,w)ds\,,
$$
with $\la\in\R$ being a parameter. We assume that $\dot{u}_\la\in
L_a^2([0,1]\times W,dt\times d\mu)$, where the subscript   ``$_a$'' means
that it is adapted to the canonical filtration for almost all $s\in
[0,1]$. We denote the primitive of $\dot{u}_\la$ by $u_\la$ and assume
that $E[\rho(-\delta u_\la)]=1$, where $\rho$ denotes the Girsanov
exponential:
$$
\rho(-\delta u_\la)=\exp\left(-\int_0^1 \dot{u}_\la(s)
  dW_s-\frac{1}{2}\int_0^1|\dot{u}_\la(s)|^2ds\right)\,.
$$
We shall assume that the map $\la\to \dot{u}_\la$ is differentiable as
a map in $L_a^2([0,1]\times W,dt\times d\mu)$, we denote its
derivative w.r. to $\la$ by $\dot{u}'_\la(s)$ or by $\dot{u}'(\la,s)$
and its primitive w.r. to $s$ is denoted as $u_\la'(t)$.
\begin{theorem}
\label{thm-1}
Suppose that   $\la\to u_\la\in
L_{loc}^p(\R,d\la;\DD_{p,1}(H))$ for some $p\geq 1$, with
$E[\rho(-\delta u_\la)]=1$ for any $\la\geq 0$ and also that 
$$
E\int_0^\la(1+\rho(-\delta u_\al)) \left|E[\delta(K_\alpha  u'_\alpha)|U_\al]\right|^pd\alpha<\infty\,,
$$
where $K_\alpha=(I_H+\nabla u_\alpha)^{-1}$. Then the map
$$
\la\to L_\la=\frac{dU_\la\mu}{d\mu}
$$
is absolutely continuous and we have 
$$
L_\la(w)=L_0\exp \int_0^\la E\Big[\delta(K_\alpha
u'_\alpha)|U_\alpha=w\Big]d\alpha\,.
$$
\end{theorem}
\nproof 
Let us note first that the map $(\la,w)\to L_\la(w)$ is measurable
thanks to the Radon-Nikodym theorem. Besides, for any (smooth) cylindrical function $f$, we have 
\beaa
\frac{d}{d\la}E[f\circ U_\la]&=&E[(\nabla f\circ U_\la,u'_\la)_H]\\
&=&E[((I_H+\nabla u_\la)^{-1\star}\nabla (f\circ U_\la),u'_\la)_H]\\
&=&E[(\nabla (f\circ U_\la),(I_H+\nabla u_\la)^{-1} u'_\la)_H]\\
&=&E[f\circ U_\la\,\delta\{(I_H+\nabla u_\la)^{-1} u'_\la\}]\\
&=&E[f\circ U_\la\,E[\delta (K_\la u'_\la)|U_\la]]\\
&=&E[f\,E[\delta (K_\la u'_\la)|U_\la=w]L_\la]\,.
\eeaa
Hence, for any fixed  $f$, we get
$$
\frac{d}{d\la}\langle f, L_\la\rangle =\langle f,L_\la E[\delta (K_\la
u'_\la)|U_\la=w]\rangle\,,
$$
both sides of the above equality  are continuous w.r. to $\la$, hence
we get
$$
<f,L_\la>-<f,L_0>=\int_0^\la <f,L_\al
E\left[\delta(K_\al u'_\al)|U_\al=w\right]>d\al\,.
$$
From the hypothesis, we have
$$
E\int_0^\la L_\al |E[\delta(K_\al
  u'_\al)|U_\al=w]|d\al=E\int_0^\la |E[\delta(K_\al
  u'_\al)|U_\al]|\,d\al<\infty\,.
$$
By the measurability of the disintegrations, the mapping  $(\al,w)\to E[\delta(K_\al
  u'_\al)|U_\al=w]$ has a measurable modification, hence  the following
  integral equation holds in the ordinary sense for almost all $w\in W$
$$
L_\la=L_0+\int_0^\la L_\al E[\delta(K_\al u'_\al)|U_\al=w] d\al\,,
$$
for $\la>0$.  Therefore the map $\la\to L_\la$ is almost surely
absolutely continuous w.r. to the Lebesgue measure. To show its
representation as an exponential, we need to show that the map 
$\al\to E[\delta(K_\al u'_\al)|U_\al =w]$ is almost surely locally
integrable. To achieve this it suffices to observe that
\beaa
E\int_0^\la | E[\delta(K_\al u'_\al)|U_\al =w]|d\al&=&E\int_0^\la |
E[\delta(K_\al u'_\al)|U_\al =w]|\frac{L_\al}{L_\al}d\al \\
&=&E\int_0^\la |
E[\delta(K_\al u'_\al)|U_\al]|\frac{1}{L_\al\circ U_\al}d\al\\
&=&E\int_0^\la |
E[\delta(K_\al u'_\al)|U_\al]|E[\rho(-\delta u_\al)|U_\al]d\al<\infty
\eeaa
by hypothesis and by Theorem \ref{thm-00}. Consequently we have
the explicit expression for $L_\la$  given as:
$$
L_\la(w)=L_0\exp \int_0^\la E[\delta(K_\al u'_\al)|U_\al=w]d\al\,.
$$
\nqed

\begin{remarkk}
\noindent
An important tool to control the hypothesis of Theorem \ref{thm-1} is  the inequality of T. Carleman which says
that (cf.  \cite{D-S}, Corollary XI.6.28)
$$
\|\dett(I_H+A)(I_H+A)^{-1}\|\leq \exp\frac{1}{2}\left(\|A\|_2^2+1\right)\,,
$$
for any Hilbert-Schmidt operator $A$, where the left hand side is the
operator norm, $\dett(I_H+A)$ denotes the modified Carleman-Fredholm
determinant and $\|\cdot\|_2$ denotes the Hilbert-Schmidt norm. Let us
remark that if $A$ is a quasi-nilpotent operator, i.e., if the
spectrum of $A$ consists of zero only, then $\dett(I_H+A)=1$, hence in
this case the Carleman inequality reads 
$$
\|(I_H+A)^{-1}\|\leq
\exp\frac{1}{2}\left(\|A\|_2^2+1\right)\,.
$$
This case happens when $A$ is equal to the Sobolev derivative of some
$u\in \DD_{p,1}(H)$ whose  drift $\dot{u}$ is
adapted to the filtration $(\calF_t,\,t\in  [0,1])$,
\end{remarkk}

 \noindent
From now on, for the sake of technical simplicity we shall assume that  {\bf{$u_\la$ is essentially bounded
 uniformly w.r.to $\la$}.}

\begin{proposition}
\label{prop-1}
Let $F\in L^p(\mu)$  then
the map $\la\to E[F|U_\la=w]$ is weakly continuous  with values in
$L^{p-}(\mu)${\footnote{$p-$ denotes any $p'<p$ and $q+$ any $q'>q$}}.
\end{proposition}
\nproof
First we have 
\beaa
\int_W|E[F|U_\la=w]|^{p}d\mu&=&\int_W|E[F|U_\la=w]|^{p}\frac{L_\la}{L_\la}d\mu\\
&=&\int_W |E[F|U_\la]|^{p-}\frac{1}{L_\la\circ U_\la}d\mu\\
&=&\int_W |E[F|U_\la]|^{p}E[\rho(-\delta u_\la)|U_\la]d\mu<\infty\,, 
\eeaa
hence $E[F|U_\la=w]\in L^{p-}(\mu)$ for any $F\in L^p(\mu)$. Besides,
for any $f\in C_b(W)$, 
$$
E[f\circ U_\la\,F]=E[f E[F|U_\la=w]\,L_\la]
$$
therefore
$$
|E[f\circ U_\la\,F]|\leq \|F\|_p \|f\circ U_\la\|_q\leq C_q
\|F\|_p\|\|f\|_{q+}\,.
$$
This relation, combined with the continuity of $\la\to f\circ U_\la$,
due to the Lusin theorem,  in $L^q$ for any $f\in L^{q+}$, implies the weak continuity of the map $\la\to
[F|U_\la=w]\,L_\la$ with  values in $L^{p-}(\mu)$, since $\la\to L_\la$
and $\la\to (L_\la)^{-1}$ are  almost surely and  strongly continuous
in $L^p(\mu)$, the claim follows.
\nqed

\begin{theorem}
\label{thm-2}
Assume that $F\in\DD_{p,1}$ for some $p>1$ and that 
$$
E\int_0^\la |\delta(FK_\al u'_\al)|d\al<\infty
$$
for any $\la>0$, then $\la\to E[F|U_\la=w]$ is $\mu$-a.s. absolutely
continuous w.r. to the Lebesgue measure $d\la$ and the map
$\la\to E[F|U_\la]$ is almost surely and hence $L^p$-continuous.
\end{theorem}
\nproof
Using the same method as in the proof of Theorem \ref{thm-1}, we
obtain
\beaa
\frac{d}{d\la}E[\theta\circ
U_\la\,F]&=&\frac{d}{d\la}E[\theta\,E[F|U_\la=w]\,L_\la]\\
&=&E[\theta \,L_\la\,E[\delta(F\,K_\la u'_\la)|U_\la=w]]
\eeaa
for any cylindrical function $\theta$. By continuity w.r.to $\la$, we
get 
$$
E\left[\theta\Big(L_\la E[F|U_\la=w]-L_0
E[F|U_0=w]\Big)\right]=\int_0^\la E\left[\theta L_\al E[\delta(FK_\al
u'_\al)|U_\al=w]\right]d\al\,.
$$
By the hypothesis
$$
E\int_0^\la | L_\al E[\delta(FK_\al u'_\al)|U_\al=w]|d\al<\infty
$$
and since $\theta$ is an arbitrary cylindrical function, we obtain the
identity 
$$
L_\la E[F|U_\la=w]-L_0 E[F|U_0=w]=\int_0^\la L_\al \,E[\delta(F K_\al
u'_\al)|U_\al=w]d\al
$$
almost surely and this proves the first part of the theorem since
$\la\to L_\la$ is already absolutely continuous and strictly
positive. For the second part, we denote $E[F|U_\la]$ by
$\hat{F}(\la)$ and we assume that $(\la_n,n\geq 1)$ tends to some
$\la$, then there exists a sub-sequence $(\hat{F}(\la_{k_l}),l\geq 1)$
which converges weakly to some limit; but, from the first part of the proof, we
know that $(E[F|U_{\la_{k_l}}=w],l\geq 1)$ converges
almost surely  to $E[F|U_\la=w]$ and by the uniform integrability,
there is also strong convergence  in $L^{p-}(\mu)$. Hence, for any cylindrical
function $G$, we have 
\beaa
E[\hat{F}(\la_{k_l})\,G]&=&E[E[F|U_{\la_{k_l}}=w]E[G|U_{\la_{k_l}}=w]L_{\la_{k_l}}]\\
&\rightarrow&E[E[F|U_{\la}=w]E[G|U_{\la}=w]L_{\la}]\\
&=&E[\hat{F}(\la)\,G]\,.
\eeaa
Consequently, the map $\la\to \hat{F}(\la)$ is weakly continuous in
$L^p$, therefore it is also strongly continuous.
\nqed

\remark
Another proof consists of remarking that 
$$
E[F|U_\la=w]|_{w=U_\la}=E[F|U_\la]
$$
$\mu$-a.s. and that $\la\to E[F|U_\la=w]$ is continuous a.s. and in
$L^{p-}$ from the first part of the proof and that $(L_\la, \la\in [a,b])$ is uniformly
integrable. These observations, combined with the Lusin's
theorem imply the continuity in $L^0(\mu)$ (i.e., in  probability)
of $\la\to E[F|U_\la]$ and the $L^p$-continuity follows.

\vspace{0.75cm}
\noindent
We shall need some technical results, to begin with, let $U^\tau_\la$
denote the shift defined on $W$ by
$$
U^\tau_\la(w)=w+\int_0^{\cdot\wedge\tau}\dot{u}_\la(s)ds\,,
$$
for $\tau\in [0,1]$. We shall denote by $L_\la(\tau)$ the
Radon-Nikodym density 
$$
\frac{dU^\tau_\la\mu}{d\mu}=L_\la(\tau)\,.
$$

\begin{lemma}
\label{lemma-1}
We have the relation
$$
L_\la(\tau)=E[L_\la|\calF_\tau]
$$
almost surely.
\end{lemma}
\nproof
Let $f$ be an $\calF_\tau$-measurable, positive, cylindrical function;
then it is straightforward to see that $f\circ U_\la=f\circ
U_\la^\tau$, hence
$$
E[f\,L_\la]=E[f\circ U_\la]=E[f\circ U_\la^\tau]=E[f\,L_\la(\tau)]\,.
$$
\nqed

\begin{lemma}
\label{lemma-2}
Let $\calU_\la^\tau(t)$ be the sigma algebra generated by
$\{U^\tau_\la(s);\,s\leq t\}$. Then, we have 
$$
E[f|\calU_\la^\tau(1)]=E[f|U_\la^\tau]
$$
for any positive, measurable function on $W$.
\end{lemma}
\nproof
Here, of course the second conditional expectation is to be understood
w.r. to the sigma algebra generated by the mapping $U_\la^\tau$ and
once this point is fixed the claim is trivial.
\nqed

\begin{proposition}
\label{prop-2}
With the notations explained above, we have
$$
L_\la(\tau)=L_0(\tau)\exp\int_0^\la E[\delta\{(I_H+\nabla
u_\al^\tau)^{-1}u'^{\tau}_\al\}|U_\al^\tau=w]d\al\,.
$$
Moreover, the map $(\la,\tau)\to L_\la(\tau)$ is continuous on
$\R\times [0,1]$ with values in $L^p(\mu)$ for any $p\geq 1$.
\end{proposition}
\nproof
The first claim can be proved as we have done in the first part of the
proof of Theorem \ref{thm-1}. For the second part, let $f$ be a positive, measurable function on W; we have
$$
E[f\circ U_\la^\tau]=E[f\,L_\la(\tau)]\,.
$$
If $(\tau_n,\la_n)\to(\tau,\la)$, from the Lusin theorem and the
uniform integrability of the densities $(L_{\la_n}(\tau_n),n\geq 1)$,
the sequence  $(f\circ U^{\tau_n}_{\la_n},\,n\geq 1)$ converges in
probability to $f\circ U_\la^\tau$, hence, again by the uniform
integrability, for any $q>1$ and $f\in L^q(\mu)$, 
$$
\lim_nE[f\,L_{\la_n}(\tau_n)]=E[f\,L_\la(\tau)]\,.
$$
From Lemma \ref{lemma-1}, we have
\beaa
E[L_{\la_n}(\tau_n)^2]&=&E[L_{\la_n}(\tau_n)\,E[L_{\la_n}|\calF_{\tau_n}]]\\
&=&E[L_{\la_n}(\tau_n)\,L_{\la_n}]\,,
\eeaa
since, from Theorem \ref{thm-1}, $L_{\la_n}\to L_\la$ strongly in all
$L^p$-spaces, it follows that $(\la,\tau)\to L_\la(\tau)$ is
$L^2$-continuous, hence also $L^p$-continuous for any $p>1$.
\nqed

\begin{proposition}
\label{prop-3}
The mapping $(\la,\tau)\to L_\la(\tau)$ is a.s. continuous, moreover
the map
$$
(\tau,w)\to (\la\to L_\la(\tau,w))
$$
is a $C(\R)$-valued continuous martingale and its restriction to
compact intervals (of $\la$) is uniformly integrable.
\end{proposition}
\nproof
Let us take the interval $\la\in [0,T]$, from Lemma \ref{lemma-1} we
have $L_\la(\tau)=E[L_\la|\calF°\tau]$, since $C([0,T])$ is a
separable Banach space and since we are working with the completed
Brownian filtration,  the latter equality implies an a.s. continuous,
$C([0,T])$-valued uniformly integrable martingale.
\nqed

\begin{theorem}
\label{thm-3}
Assume that 
$$
E\int_0^\la\int_0^1\left(|\delta(\dot{u}_\al(s)K_\al
  u'_\al)|+|\dot{u}'_\al(s)|^2\right)ds<\infty
$$
for any $\la\geq 0$, then the map
$$
\la\to E[\dot{u}_\la(t)|\calU_\la(t)]
$$
is continuous with values in $L_a^p(\mu,\,L^2([0,1],\R^d))$, $p\geq
1$.
\end{theorem}
\nproof
Let $\xi\in L_a^\infty(\mu,H)$ be  smooth and cylindrical, then, by
similar calculations as in the proof of Theorem \ref{thm-2}, we get
\beaa
\frac{d}{d\la}E[(\xi\circ U_\la,u_\la)_H]&=&\frac{d}{d\la}<\xi\circ
U_\la,u_\la>=
\frac{d}{d\la}<\xi\circ U_\la,\hat{u}_\la>\\
&=&E\int_0^1\dot{\xi}_sL_\la(s)E\left[\delta(\dot{u}_\la(s)K_\la
u'_\la)+\dot{u}'_\la(s)|U_\la^s=w\right]ds\,,
\eeaa
but the l.h.s. is equal to
$$
E[(\nabla\xi\circ U_\la[u'_\la],u_\la)_H+(\xi\circ U_\la,u'_\la)_H]\,,
$$
which is  continuous w.r. to $\la$ provided that $\xi$ is smooth,
and that $\la\to (u'_\la, u_\la)$ is continuous in $L^p$ for $p\geq
2$. Consequently, we have the relation 
$$
<\xi\circ U_\la,u_\la>-<\xi\circ U_0,u_0>=
E\int_0^\la \int_0^1\dot{\xi}_sL_\al(s)E\left[\delta (\dot{u}_\al(s)K_\al
u'_\al)+\dot{u}'_\al(s)|U_\al^s=w\right]dsd\al
$$
and the hypothesis implies that $\la\to
L_\la(s)E[\dot{u}_\la(s)|U_\la^s=w]$ is $\mu$-a.s. absolutely continuous
w.r.to the Lebesgue measure $d\la$. Since $\la\to L_\la(s)$ is also
a.s. absolutely continuous, it follows that $\la\to
E[\dot{u}_\la(s)|U_\la^s=w]$ is a.s. absolutely continuous. Let us
denote this disintegration as the kernel $N_\la(w,\dot{u}_\la(s))$,
then
$$
N_\la(U_\la^s(w),\dot{u}_\la(s))=E[\dot{u}_\la(s)|U_\la^s]
$$
a.s. From the Lusin theorem, it follows that the map $\la\to
N_\la(U_\la^s,\dot{u}_\la(s))$ is continuous with values in
$L_a^0(\mu,L^2([0,1],\R^d))$ and the $L^p$-continuity follows from the
dominated convergence theorem.
\nqed

\begin{remarkk}
In the proof above we have the following result: assume that $\la\to
f_\la$ is continuous in $L^0(\mu)$, then $\la\to f_\la\circ U_\la$ is
also continuous in $L^0(\mu)$ provided that the family
$$
\left\{\frac{dU_\la\mu}{d\mu},\la\in [a,b]\right\}
$$
is uniformly integrable for any compact interval $[a,b]$. To see this, it suffices to verify the
sequential continuity; hence assume that $\la_n\to \la$, then we have
\beaa
\mu\{|f_{\la_n}\circ U_{\la_n}-f_{\la}\circ
U_{\la}|>c\}&\leq&\mu\{|f_{\la_n}\circ U_{\la_n}-f_{\la}\circ
U_{\la_n}|>c/2\}\\
&&+\mu\{|f_{\la}\circ U_{\la_n}-f_{\la}\circ
U_{\la_n}|>c/2\}\,,
\eeaa
but 
$$
\mu\{|f_{\la_n}\circ U_{\la_n}-f_{\la}\circ
U_{\la_n}|>c/2\}=E[L_{\la_n}1_{\{|f_{\la_n}-f_\la|>c/2\}}]\to 0
$$
by the uniform integrability of $(L_{\la_n},\,n\geq 1)$ and the
continuity of $\la\to f_\la$. The second term tends  also to zero by
the standard use of Lusin theorem and again by the  the uniform integrability of $(L_{\la_n},\,n\geq 1)$.
\end{remarkk}

\begin{corollary}
\label{cor-1}
The map $\la\to E[\rho(-\delta u_\la)|U_\la]$ is continuous as an
$L^p(\mu)$-valued map for any $p\geq 1$.
\end{corollary}
\nproof
We know that 
$$
E[[\rho(-\delta u_\la)|U_\la]=\frac{1}{L_\la\circ U_\la}\,.
$$
\nqed

\begin{corollary}
\label{cor-2}
Let $Z_\la(t)$ be the innovation process associated to $U_\la$, then 
$$
\la\to \int_0^1E[\dot{u}_\la(s)|\calU_\la(s)]dZ_\la(s)
$$
is continuous as an $L^p(\mu)$-valued map for any $p\geq 1$.
\end{corollary}
\nproof
We have 
$$
\log L_\la\circ
U_\la=\int_0^1E[\dot{u}_\la(s)|\calU_\la(s)]dZ_\la(s)+\frac{1}{2}\int_0^1|E[\dot{u}_\la(s)|\calU_\la(s)]|^2ds\,,
$$
since the l.h.s. of this equality and the second  term at the right are
continuous, the first  term at the right should be also continuous.
\nqed

\begin{theorem}
\label{thm-4}
Assume that 
$$
E\int_0^\la|\delta\{\delta(K_\al u'_\al)K_\al u'_\al-K_\al\nabla
u'_\al K_\al u'_\al+K_\al u''_\al\}|d\al<\infty
$$
for any $\la\geq 0$. Then the map 
$$
\la\to \frac{d}{d\la}L_\la
$$
is a.s. absolutely continuous w.r.to the Lebesgue measure $d\la$ and
we have 
$$
\frac{d^2}{d\la^2}L_\la(w)=L_\la E[\delta D_\la|U_\la=w]\,,
$$
where 
$$
D_\la=\delta(K_\la u'_\la)K_\la u'_\la-K_\la\nabla u'_\la K_\la
u'_\la+K_\la u''_\la\,.
$$
\end{theorem}
\nproof
Let $f$ be a smooth function on $W$, using the integration by parts
formula as before, we get
\beaa
\frac{d^2}{d\la^2}E[f\circ U_\la]&=&\frac{d}{d\la}
E[f\circ U_\la\,\delta(K_\la u'_\la)]\\
&=&E[(\nabla f\circ U_\la,u'_\la)_H\delta(K_\la u'_\la)]\\
&=&E[(K_\la^\star\nabla(f\circ U_\la),u'_\la)_H\delta(K_\la
u'_\la)+f\circ U_\la\delta(-K_\la\nabla u'_\la K_\la u'_\la+K_\la
u''_\la)]\\
&=&E\left[f\circ U_\la\left\{\delta(\delta(K_\la u'_\la)K_\la
u'_\la)-\delta(K_\la\nabla u'_\la K_\la u'_\la)+\delta(K_\la
u''_\la)\right\}\right]\,.
\eeaa
Let us define the map $D_\la$ as
$$
D_\la=\delta(K_\la u'_\la)K_\la u'_\la-K_\la\nabla u'_\la K_\la
u'_\la+K_\la u''_\la\,,
$$
we have obtained then the following relation
$$
\frac{d^2}{d\la^2}E[f\circ U_\la]=E[f\,L_\la\,E[\delta D_\la|U_\la=w]]
$$
hence
$$
< \frac{d}{d\la}L_\la,f>-<\frac{d}{d\la}L_\la,f>|_{\la=0}=\int_0^\la E[f\,L_\al\,E[\delta
D_\al|U_\al=w]]d\al\,.
$$
The hypothesis implies the existence of the strong (Bochner) integral
and we conclude that 
$$
L'_\la-L'_0=\int_0^\la L_\al E[\delta D_\al|U_\al=w]d\al
$$
a.s. for any $\la$, where $L'_\la$ denotes the derivative of $L_\la$
w.r.to $\la$.
\nqed

\begin{theorem}
\label{thm-5}
Define the sequence of functionals  inductively as
\beaa
D_\la^{(1)}&=&D_\la\\
D_\la^{(2)}&=&(\delta D_\la^{(1)})K_\la
u'_\la+\frac{d}{d\la}D_\la^{(1)}\\
&&\ldots\\
D_\la^{(n)}&=&(\delta D_\la^{(n-1)})K_\la
u'_\la+\frac{d}{d\la}D_\la^{(n-1)}\,.
\eeaa
Assume that 
$$
E\int_0^\la|\delta D^{(n)}_\al|d\al<\infty
$$
for any $n\geq 1$ and $\la\in \R$, then $\la\to L_\la$ is almost
surely a $C^\infty$-map and denoting by $L^{(n)}_\la$ its derivative
of order $n\geq 1$, we have
$$
L_\la^{(n+1)}(w)-L_0^{(n+1)}(w)=\int_0^\la L_\al  E[\delta
D_\al^{(n)}|U_\al=w]d\al\,.
$$
\end{theorem}

\section{\bf{Applications to the invertibility of adapted perturbations of
  identity}}
\label{S-Inv}
Let $u\in L_a^2(\mu, H)$, i.e., the space of square integrable,
$H$-valued functionals whose Lebesgue density, denoted as $\dot{u}(t)$, is adapted to the
filtration $(\calF_t,t\in [0,1])$ $dt$-almost surely. A frequently
asked question ire the conditions which imply the almost sure
invertibility of the adapted perturbation of identity (API) $w\to
U(w)=w+u(w)$. The next theorem gives such a condition:

\begin{theorem}
\label{thm-8}
Assume that $u\in L_a^2(\mu, H)$ with $E[\rho(-\delta u)]=1$, let
$u_\al$ be defined as $P_\al u$, where $P_\al=e^{-\al \calL}$ denotes the
Ornstein-Uhlenbeck semi-group on the Wiener space. If there exists a
$\la_0$ such that 
\beaa
\lefteqn{E\int_0^\la E[\rho(-\delta u_\al)|U_\la]\Big|E[\delta(K_\al
u'_\al)|U_\al]\Big|d\al}\\
&&=E\int_0^\la E[\rho(-\delta
u_\al)|U_\la]\Big|E[\delta((I_H+\nabla u_\al)^{-1}\calL u_\al)|U_\al]\Big|d\al<\infty
\eeaa
for $\la\leq \la_0$, then $U$ is almost surely invertible. In
particular the functional stochastic differential equation
\beaa
dV_t(w)&=&-\dot{u}(V_s(w),s\leq t)dt+dW_t\\
V_0&=&0
\eeaa
has a unique strong solution.
\end{theorem}
\nproof
Since $u_\al$ is an $H-C^\infty$-function, cf. \cite{BOOK}, the API
$U_\al=I_W+u_\al$ is a.s. invertible, cf.\cite{INV}, Corollary 1. By the hypothesis
and from Lemma 2 of \cite{INV}, $(\rho(-\delta u_\al),\al\leq \la_0)$
is uniformly integrable. Let $L_\al$ and $L$  be respectively  the Radon-Nikodym derivatives
of $U_\al\mu$ and $U\mu$ w.r. to $\mu$. From Theorem \ref{thm-1}, 
$$
L_\la(w)=L(w)\,\exp\int_0^\la E[\delta(K_\al u'_\al)|U_\al=w]d\al
$$
for any $\la\leq \la_0$ and also that $\int_0^\la |E[\delta(K_\al
u'_\al)|U_\al=w]|d\al<\infty$ almost surely. Consequently
$$
L_\la-L=\left(\exp\int_0^\la E[\delta(K_\al  u'_\al|U_\al=w]d\al-1\right) L \to 0
$$
as $\la\to 0$, in probability (even in $L^1$). We claim that the set
$(L_\al \log L_\al, \al\leq \la_0)$ is uniformly integrable. To see
this let $A\in\calF$, then 
\beaa
E[1_A L_\al \log L_\al]&=&E[1_A\circ U_\al\,\log L\circ U_\al]\\
&=&-E[1_A\circ U_\al\,\log E[\rho(-\delta u_\al)|U_\al]]\\
&\leq&-E[1_A\circ U_\al\,\log \rho(-\delta u_\al)]\\
&=&E\left[1_A\circ U_\al\left(\delta u_\al+\frac{1}{2}|u_\al|_H^2\right)\right]
\eeaa
Since $(|u_\al|^2,\al\leq \la_0)$ is uniformly integrable, for any
given $\eps>0$, there exists some $\ga>0$, such that $\sup_\al
E[1_B|u_\al|^2]\leq \eps$ as soon as $\mu(B)\leq \ga$ and this happens
uniformly w.r. to $B$, but as $(L_\al ,\al\leq \la_0)$ is uniformly
integrable, there exists a $\ga_1>0$ such that, for any $A\in \calF$,
with $\mu(A)\leq \ga_1$, we have $\mu(U_\al^{-1}(A))\leq \ga$ uniformly
in $\al$ and we obtain $E[1_A\circ U_\al |u_\al|_H^2]\leq \eps$ with
such a choice of $A$. For the first term above we have 
$$
E[1_A\circ U_\al \delta u_\al]\leq E[1_A
L_\al]^{1/2}\|u_\al\|_{L^2(\mu,H)}\leq \eps
$$
again by the same reasons. Hence we can conclude that 
$$
\lim_{\al\to 0}E[L_\al \log L_\al]=E[L \log L]\,.
$$
Moreover, as shown in \cite{ASU-2,ASU-3}, the invertibility of $U_\al$ is equivalent to 
$$
E[L_\al \log L_\al]=\frac{1}{2}E[|u_\al|_H^2]\to
\frac{1}{2}E[|u|_H^2]\,,
$$
therefore 
$$
E[L \log L]=\frac{1}{2}E[|u|_H^2]
$$
which is a necessary and sufficient condition for the invertibility of $U$

\nqed

\noindent
In several applications we encounter a situation as follows: assume
that $u:W\to H$ is a measurable map with the following property 
$$
|u(w+h)-u(w)|_H\leq c|h|_H
$$
a.s., for any $h\in H$, where $0<c<1$ is a fixed constant, or
equivalently  an upper bound like $\|\nabla u\|_{op}\leq c$ where
$\|\cdot\|_{op}$ denotes the operator norm. Combined with some
exponential integrability of the Hilbert-Schmidt norm $\nabla u$, one
can prove the invertibility of $U=I_W+u$, cf. Chapter 3 of
\cite{BOOK}. Note that  the hypothesis
$c<1$ is indispensable because of the fixed-point techniques used to
construct the inverse of $U$. However, using  the techniques developed in this paper
we can relax this rigidity of the theory:

\begin{theorem}
\label{thm-9}
Let $U_\la=I_W+\la u$ be an API (adapted perturbation of identity)
with $u\in \DD_{p,1}(H)\cap L^2(\mu,H)$, such that, for any $\la<1$, $U_\la$ is
a.s. invertible. Assume that 
\begin{equation}
\label{con-suff}
E\int_0^1\rho(-\delta (\al u))|E[\delta((I_H+\al\nabla
u)^{-1}u)|U_\al]|d\al<\infty\,.
\end{equation}
Then $U=U_1$ is also a.s. invertible.
\end{theorem}
\nproof
Let $L=L_1$ be the Radon-Nikodym derivative of $U_1\mu$ w.r. to
$\mu$. It suffices to show that 
$$
E[L\log L]=\frac{1}{2}E[|u|_H^2]
$$
which is an equivalent condition to the a.s.  invertibility of $U$,
cf. \cite{ASU-3}. For this it suffices to show first that
$(L_\la,\la<1)$ converges in $L^0(\mu)$ to $L$, then that $(L_\la \log
L_\la,\la<1)$ is uniformly integrable. The first claim follows from
the hypothesis (\ref{con-suff}) and the second claim can be proved
exactly as in the proof of Theorem \ref{thm-8}.
\nqed

\section{\bf{Variational applications to  entropy and estimation}}
\label{S-Ent}
In the estimation and information  theories, one  often encounters the problem of
estimating the signal $u_\la$ from the observation data generated by
$U_\la$ and then verifies the various properties of the mean square
error w.r.to the signal to noise ratio, which is represented in our
case with the parameter $\la$. Since we know that (\cite{ASU-3})
$$
E[L_\la\log L_\la]=\frac{1}{2}E\int_0^1|E[\dot{u}_\la(s)|\calU_\la(s)]|^2ds\,,
$$
the behavior of the  mean square error is completely characterized by
that of the relative entropy. Let $\theta$ denote the entropy of
$L_\la$ as a function of $\la$:
$$
\theta(\la)=E[L_\la\log L_\la]\,.
$$
From our results, it comes immediately
that 
\beaa
\frac{d\theta(\la)}{d\la}&=&E[L'_\la\log L_\la]\\
&=&E[L_\la\,E[\delta(K_\la u'_\la)|U_\la=w]\log L_\la]\\
&=&E[E[\delta(K_\la u'_\la)|U_\la]\log L_\la\circ U_\la]\\
&=&-E[\delta(K_\la u'_\la)\log E[\rho(-\delta u_\la)|U_\la]]\,.
\eeaa
Similarly
\beaa
\frac{d^2\theta(\la)}{d\la^2}&=&E\left[L''_\la\log L_\la
+(L'_\la)^2\frac{1}{L_\la}\right]\\
&=&E[L''_\la\log L_\la+L_\la\,E[\delta(K_\la u'_\la)|U_\la=w]^2]\\
&=&E[E[\delta D_\la|U_\la=w]L_\la\log L_\la++L_\la\,E[\delta(K_\la
u'_\la)|U_\la=w]^2]\\
&=&E[E[\delta D_\la|U_\la]\log L_\la\circ U_\la+E[\delta(K_\la u'_\la)|U_\la]^2]\,.
\eeaa

\noindent
In particular we have
\begin{theorem}
\label{thm-6}
Assume that 
$$
E\left[E[\delta
D_\la|U_\la]\left(\int_0^1E[\dot{u}_\la(s)|\calU_\la(s)]dZ_\la(s)+\frac{1}{2}\int_0^1|E[\dot{u}_\la(s)|\calU_\la(s)]|^2ds\right)\right]<
E\left[E[\delta(K_\la u'_\la)|U_\la]^2\right]
$$
for some $\la=\la_0>0$, then there exists an $\eps>0$ such that the {\bf entropy is convex} as a function of $\la$ on the
interval $(\la_0-\eps,\la_0+\eps)$. In particular, if $u_0=0$, then
the same conclusion holds true on some $(0,\eps)$.
\end{theorem}
% \nproof
% Note that $D_0=u\delta u-\nabla u[u]$, hence $\delta D_0=(\delta
% u)^2-|u|_H^2$ and we get 
% $$
% \theta''(0)=E[(\delta u)^2]>0\,.
% $$
% Let 
% $$
% \eps=\inf(\la>0:\,\theta''(\la)\leq 0)\,,
% $$
% then the proof follows by continuity.
% \nqed
% \

\subsection{Applications to the anticipative estimation}
In this section we study briefly the estimation of $\dot{u}_\la(t)$
with respect to the final filtration $\calU_\la(1)=\sigma(U_\la)$.
\begin{theorem}
\label{thm-7}
Assume that 
$$
E\int_0^\la L_\al |E[\dot{u}'_\al(s)+\delta(\dot{u}_\al(s)K_\al
u'_\al)|U_\al]|^p d\al<\infty\,,
$$
for a $p\geq 1$, then, $dt$-a.s., the map $\la\to L_\la E[\dot{u}_\la(t)|U_\la=x]$ and
hence the map  $\la\to
E[\dot{u}_\la(t)|U_\la=x]$ are strongly differentiable in $L^p(\mu)$
for any $p\geq 1$ and we have 
$$
\frac{d}{d\la}E[\dot{u}_\la
(t)|U_\la=x]=E[\dot{u}'_\la(t)+
\delta(\dot{u}_\la(t)K_\la u'_\la)|U_\la=x]-E[\dot{u}_\la(t)|U_\la=x] E[\delta(K_\la u'_\la)|U_\la=x]
$$
$d\mu\times dt$-a.s.
\end{theorem}
\nproof
For a smooth function $h$ on $W$, we have
\beaa
\frac{d}{d\la}<E[\dot{u}_\la
(t)|U_\la=x],h\,L_\la>&=&\frac{d}{d\la}<E[\dot{u}_\la
(t)|U_\la],h\circ U_\la>\\
&=&E[\dot{u}'_\la(t)h\circ U_\la+\dot{u}_\la(t)(\nabla g\circ
U_\la,u'_\la)_H]\\
&=&E[E[\dot{u}'_\la(t)|U_\la]h\circ U_\la+h\circ
U_\la\delta(\dot{u}_\la(t)K_\la u'_\la)]\\
&=&E\left[hL_\la(x)\left(E[\dot{u}'_\la(t)|U_\la=x]+E[\delta(\dot{u}_\la(t)K_\la
u'_\la)|U_\la=x]\right)\right]\,.
\eeaa
The hypothesis implies that this weak derivative is in fact a strong
one in $L^p(\mu)$, the formula follows by dividing both sides by
$L_\la$ and by the explicit form of $L_\la$ given in Theorem
\ref{thm-1}.
\nqed

Using the formula of Theorem \ref{thm-7}, we can study the behavior
of the error of  non-causal estimation of $u_\la$ (denoted as NCE in
the sequel)  defined as 
\beaa
NCE&=&E\int_0^1|\dot{u}_\la(s)-E[\dot{u}_\la(s)|\calU_\la(1)]|^2ds\\
&=&E\int_0^1|\dot{u}_\la(s)-E[\dot{u}_\la(s)|U_\la]|^2ds
\eeaa
To do this we prove some technical results:
\begin{lemma}
\label{lem-1}
Assume that 
\begin{equation}
\label{hyp-1}
E\int_0^\la\int_0^1|\dot{u}''_\al(s)+\delta(\dot{u}'_\al(s)K_\al u'_\al)|^pds
d\al<\infty
\end{equation}
for some $p>1$, for any $\la>0$, then the map 
$$
\la\to L_\la E[\dot{u}'_\la(s)|U_\la=x]
$$
is strongly differentiable in $L_a^p( d\mu,L^2([0,1]))$, and its
derivative is equal to
$$
L_\la E[\dot{u}''_\la(s)+\delta (\dot{u}_\la'(s)K_\la
u'_\la)|U_\la=x]
$$
$ds\times d\mu$-a.s.
\end{lemma}
\nproof
Let $h$ be a cylindrical function on $W$, then, using, as before, the
integration by parts formula, we get
\beaa
\frac{d}{d\la}E[L_\la
E[\dot{u}'_\la(s)|U_\la=x]\,h]&=&\frac{d}{d\la}E[\dot{u}_\la'(s)\,h\circ
U_\la]\\
&=&E[\dot{u}''_\la(s)\,h\circ U_\la+h\circ
U_\la\,\delta(\dot{u}_\la'(s)K_\la u'_\la)]\\
&=&E\left[h\,L_\la\left(E[\dot{u}''_\la(s)+\delta (\dot{u}_\la'(s)K_\la
u'_\la)|U_\la=x]\right)\right]\,.
\eeaa
This proves that the weak derivative satisfies the claim, the fact
that it coincides with the strong derivative follows from the
hypothesis (\ref{hyp-1}).
\nqed

\noindent
Let us define the variance  of the estimation as  
$$
\beta(\la,s)=E\left[|E[\dot{u}_\la(s)|\calU_\la(1)]|^2\right]\,,
$$
we shall calculate the first two derivatives of $\la\to \beta(\la,s)$ w.r.to
$\la$ in order  to observe its variations. Using Lemma \ref{lem-1}, we
have immediately the first derivative as
\bea
\frac{d}{d\la}\beta(\la,s)&=&E\Bigg[ E[\dot{u}_\la(s)|U_\la=x]L_\la \nonumber \\
&& \left(E[\dot{u}'_\la(s)+\delta(\dot{u}_\la(s)K_\la
u'_\la)|U_\la=x]-\frac{1}{2}E[\dot{u}_\la(s)|U_\la=x]E[\delta(K_\la
u'_\la)|U_\la=x]\right)\Bigg]
\label{der-1}
\eea

\noindent
The proof of the following lemma can be done exactly in the same
manner as before, namely, by verifying first the weak
differentaibility using cylindrical functions and then assuring that
the hypothesis implies the existence of the strong derivative and it
is left to the reader:

\begin{lemma}
\label{lem-2}
Assume that 
$$
E\int_0^\la |\delta(\delta(K_\al u'_\al)K_\al u'_\al)+\delta(K_\al
u''_\al-K_\al \nabla u'_\al K_\al u'_\al)|^pd\al<\infty\,,
$$
for some $p\geq 1$.
Then the map
$$
\la\to L_\la E[\delta(K_\la u'_\la)|U_\la=x]
$$
is strongly differentiable in $L^p(\mu)$ and we have 
\beaa
\frac{d}{d\la}(L_\la E[\delta(K_\la u'_\la)|U_\la=x])&=&L_\la
E\left[\delta(\delta(K_\la u'_\la)K_\la u'_\la)|U_\la=x\right]\\
&&+L_\la E\left[\delta(K_\la u''_\la -K_\la\nabla u'_\la K_\la
u'_\la)|U_\la=x\right]\,.
\eeaa
\end{lemma}

\noindent
Combining Lemma \ref{lem-1} and Lemma \ref{lem-2} and including the
action of $L_\la$, we conclude that
\beaa
\beta''(\la)&=&E\Big[E[\dot{u}''_\la +\delta(\dot{u}'_\la K_\la
u'_\la)|U_\la]E[\dot{u}_\la(s)|U_\la]\Big]\\
&&+E\Big[E[\dot{u}'_\la(s)|U_\la]\Big(E[\dot{u}'_\la(s)+\delta(\dot{u}_\la(s)K_\la
  u'_\la)|U_\la]\\
&&\,\,\,\,\,\,-E[\dot{u}_\la(s)|U_\la]E[\delta(K_\la
u'_\la)|U_\la]\Big)\Big]\\
&&+E\Big[E[\delta\left\{\dot{u}''_\la(s)K_\la
u'_\la-\dot{u}'_\la(s)K_\la \nabla u'_\la K_\la
u'_\la\right\}\\
&&\,\,\,\,\,\,+\delta\left\{\dot{u}_\la (s)K_\la u''_\la+\delta(\dot{u}_\la
(s)K_\la u'\la)K_\la
u'_\la\right\}|U_\la]E[\dot{u}_\la(s)|U_\la]\Big]\\
&&+E\Big[E[\delta(\dot{u}_\la(s)K_\la
  u'_\la)|U_\la]\Big(E[\dot{u}'_\la(s) +\delta(\dot{u}_\la(s) K_\la
    u'_\la)|U_\la]\\
&&\,\,\,\,\,\,-E[\dot{u}_\la(s)|U_\la]E[\delta(K_\la
    u'_\la]|U_\la]\Big)\Big]\\
&&-E\Big[E[\dot{u}_\la(s)|U_\la]\left(E[\dot{u}'_\la(s)+\delta(\dot{u}_\la(s)
  K_\la u'_\la)|U_\la]-E[\dot{u}_\la(s)|U_\la] E[\delta(K_\la
  u'_\la)|U_\la]\right)E[\delta(K_\la u'_\la)|U_\la]\Big]\\
&&\,\,\,\,\,\,-\frac{1}{2}E\Big[E[E[\dot{u}_\la
  (s) U_\la]^2\left\{E[\delta(\delta(K_\la u'_\la)K_\la u'_\la+K_\la
  u''_\la-K_\la\nabla u'_\la K_\la u'_\la)|U_\la]\right\}\Big]\,.
\eeaa

Assume now that $\la\to u_\la$ is linear, then a simple calculation
shows  that 
$$
\beta''(0)=E[|\dot{u}(s)|^2]\,,
$$
hence the quadratic norm of the non-causal estimation of $u$, i.e., the function
$$
\la\to E\int_0^1|E[\dot{u}_\la (s)|\calU_\la(1)]|^2ds
$$
is convex at some vicinity of $\la=0$.

\subsection{Relations with Monge-Kantorovich measure transportation}
Since $L_\la\log L_\la\in L^1(\mu)$, it follows the existence of
$\phi_\la\in \DD_{2,1}$, which is $1$-convex (cf. \cite{F-U0}) such that $(I_W+\nabla
\phi_\la)\mu=L_\la\cdot \mu$ (i.e., the measure with density
$L_\la$), cf. \cite{F-U1}. From the $L^p$-continuity of the map $\la\to L_\la$ and from the dual
characterization of the Monge-Kantorovich problem, \cite{Vil}, we
deduce the measurability of the transport potential $\phi_\la$ as a
mapping of $\la$. Moreover there exists a non-causal Girsanov-like
density $\La_\la$ such that 
\begin{equation}
\label{maq}
\La_\la\,L_\la\circ T_\la=1
\end{equation}
$\mu$-a.s., where $\La_\la$ can be expressed as 
$$
\La_\la=J(T_\la)\exp\left(-\frac{1}{2}|\nabla \phi_\la|_H^2\right)\,,
$$
where $T_\la\to J(T_\la)$ is a log-concave, normalized determinant
(cf.\cite{F-U2}) with values in $[0,1]$. Using the relation (\ref{maq}), we obtain another
expression for the entropy:
\beaa
E[L_\la\log L_\la]&=&E[\log L_\la\circ T_\la]\\
&=&-E[\log \La_\la]\\
&=&E\left[-\log J(T_\la)+\frac{1}{2}|\nabla \phi_\la|_H^2\right]\,.
\eeaa
Consequently, we have 
\beaa
\frac{1}{2}E\int_0^1|E[\dot{u}_\la(s)\mid
\calU_\la(s)]|^2ds&=&E\left[-\log J(T_\la)+\frac{1}{2}|\nabla
  \phi_\la|_H^2\right]\\
&=&E\left[-\log J(T_\la)\right]+\frac{1}{2}d^2_H(\mu,L_\la\cdot\mu)\,,
\eeaa
where $d_H(\mu,L_\la\cdot\mu)$ denotes the Wasserstein distance along
the Cameron-Martin space between the  probability measures $\mu$ and
$L_\la\cdot \mu$. This result gives another explanation for the
property remarked in \cite{M-Z} about the independence of the
quadratic norm of the estimation from the filtrations with respect to
which the causality notion  is defined. Let us remark finally that if 
$$
d_H(\mu,L_\la\cdot\mu)=0
$$
then $L_\la=1$ $\mu$-almost surely hence $E[\dot{u}_\la(s)\mid
\calU_\la(s)]=0$ $ds\times d\mu$-a.s. Let us note that such a case
may happen without having $u_\la=0$ $\mu$-a.s. As an example let us
choose an API, say $K_\la=I_W+k_\la$ which is not almost surely
invertible for any $\la\in (0,1]$. Assume that $E[\rho(-\delta
k_\la)]=1$ for any $\la$. We have 
$$
\frac{dK_\la\mu}{d\mu}=\rho(-\delta m_\la)
$$
for some $m_\la\in L^0_a(\mu,H)$, define $M_\la=I_W+m_\la$, then
$U_\la=M_\la\circ K_\la$ is a Brownian motion and an API, hence (cf. \cite{ASU-4}) it
should be equal to its own innovation process and this is equivalent
to say that $E[\dot{u}_\la(s)\mid \calU_\la(s)]=0$ $ds\times d\mu$-a.s.
\section{\bf{Applications to Information Theory}}
\label{S-Info}
\noindent
In this section we give first an extension of the results about the
quadratic error in the additive nonlinear Gaussian model which extends
the results of \cite{D,G-Y,KZZ,M-Z} in the sense that we drop a basic
assumption made implicitly or explicitly in these works; namely  the conditional
form of the signal is not   an invertible perturbation of
identity. Afterwards we study the variation of this quadratic error
with respect to a parameter on whose depends the information channel
in a reasonably smooth manner.

Throughout  this section we shall suppose  the existence of the signal in
the following form:
$$
U(w,m)=w+u(w,m)
$$
where $m$ runs in a measurable space $(M,\calM)$ governed with a
measure $\nu$ and independent of the Wiener path $w$, later on we
shall assume that the above signal  is also parametrized with a scalar $\la\in \R$. We suppose also
that, for each fixed $m$, $w\to U(w,m)$ is an adapted perturbation of
identity with $E_\mu[\rho(-\delta u(\cdot,m))]=1$ and that 
$$
\int_0^1 \int_{W\times M}|\dot{u}_s(w,m)|^2ds d\nu d\mu<\infty\,.
$$
In the sequel we shall denote the product measure $\mu\otimes \nu$ by
$\ga$ and $P$ will represent the image of $\ga$ under the map
$(w,m)\to (U(w,m),m)$, moreover we shall denote by $P_U$ the first
marginal of $P$.

The following result is known in several different cases,
cf. \cite{D,G-Y,KZZ,M-Z},  and we give
its proof in the most general case:

\begin{theorem}
\label{info-1}
Under the assumptions explained above the following relation between
the mutual information $I(U,m)$ and the quadratic estimation error  holds
true:
$$
I(U,m)=\int_{W\times M}\log\frac{dP}{dP_U\otimes d\nu} dP=
\frac{1}{2}E_\ga\int_0^1\Big(|E_\mu[\dot{u}_s(w,m)|\calU_s(m)]|^2-|E_\ga[\dot{u}_s|\calU_s]|^2\Big)ds\,,
$$
where $(\calU_s(m),\,s\in [0,1])$ is the filtration  generated by the partial map
$w\to U(w,m)$.
\end{theorem}
\proof
Let us note that the map $(s,w,m)\to E_\mu[f_s|\calU_s(m)]$ is measurable for
any positive, optional  $f$. To proceed to the proof,  remark
first that 
\bea
\frac{dP}{dP_U\otimes d\nu}&=&\frac{dP}{d\ga}\,\frac{d\ga}{dP_U\otimes d\nu}      \label{rel-1}\\
\frac{d\ga}{dP_U\otimes d\nu}&=&\frac{d\mu\otimes d\nu}{dP_U\otimes
  d\nu}=\left(\frac{dP_U}{d\mu}\right)^{-1} \label{rel-2}
\eea
since $P_U\sim \mu$. Think of $w\to U(w,m)$ as an API on the Wiener
space for each fixed $m\in M$. The image of the Wiener measure $\mu$
under  this map is absolutely continuous w.r. to $\mu$; denote the
corresponding density as $L(w,m)$. We have for any positive,
measurable function $f$ on $W\times M$
\beaa
E_P[f]&=&E_\ga[f\circ U]\\
&=&\int_{W\times M }f( U(w,m),m) d\nu(m) d\mu(w)\\
&=&\int_M E_\mu\left[f\frac{dU(\cdot,m)\mu}{d\mu}\right]d\nu(m)\\
&=&E_\ga[f L]\,,
\eeaa
hence $(w,m)\to L(w,m)$ is the Radon-Nikodym density of $P$ w.r. to
$\ga$. From \cite{ASU-3} we have at once
$$
E_\mu[L(\cdot,m)\log
L(\cdot,m)]=\frac{1}{2}E_\mu\int_0^1|E_\mu[\dot{u}_s(\cdot,m)|\calU_s(m)]|^2ds\,.
$$
Calculation of  $dP_U/d\mu$ is immediate: 
$$
\hat{L}=\frac{dP_U}{d\mu}(w)=\int_M L(w,m)d\nu(m).
$$
Moreover from the Girsanov theorem, we have
$$
E_\ga[f\circ U\,\rho(-\delta u(\cdot,m))]=E_\ga[f]
$$
for any $f\in C_b(W)$. Denote by $\calU_t$ the sigma algebra generated
by $(U_s:\,s\leq t)$ on $W\times M$. It is easy to see that the
process $Z=(Z_t,t\in [0,1])$, defined by 
$$
Z_t=U_t(w,m)-\int_0^t E_\ga[\dot{u}_s|\calU_s]ds
$$
is a $\ga$-Brownian motion and any $(\calU_t,\,t\in [0,1])$- local
martingale w.r. to $\ga$  can be
represented as a stochastic integral w.r. to the innovation process
$Z$, cf. \cite{FKK}. Let $\hat{\rho}$ denote 
\begin{equation}
\label{ro}
\hat{\rho}=\exp\left(-\int_0^1E_\ga[\dot{u}_s|\calU_s]dZ_s-\frac{1}{2}\int_0^1|E_\ga[\dot{u}_s|\calU_s]|^2ds\right)
\end{equation}
Using again the Girsanov theorem we obtain the following equality 
$$
E_\ga\left[f\circ U \hat{\rho}\right]=E_\ga\left[f\circ U \rho(-\delta u(w,m))\right]
$$
for any nice $f$.  This result  implies that
$$
E_\ga[\rho(-\delta u)|U]=\hat{\rho}
$$
$\ga$-almost surely. Besides, for nice $f$ on $W$, 
\beaa
E_{P_U}[f]&=&E_\ga[f\circ U]=E_\ga[f L]=E_\ga[f \hat{L}]\\
&=&E_\ga[f\circ U\hat{L}\circ U\,\rho(-\delta u)]\\
&=&E_\ga[f\circ U\hat{L}\circ U\,\hat{\rho}]
\eeaa
which implies that
$$
\hat{L}\circ U\,\hat{\rho}=1
$$
$\ga$-almost surely. We have calculated all the necessary ingredients
to prove the claimed representation of  the mutual
information $I(U,m)$:
\beaa
I(U,m)&=&E_P\left[\log\left(\frac{dP}{d\ga}\cdot \frac{d\ga}{dP_U\otimes d\nu}\right)\right]\\
&=&E_P\left[\log \frac{dP}{d\ga}+\log \frac{d\ga}{dP_U\otimes d\nu}\right]\\
&=&E_\ga\left[\frac{dP}{d\ga} \log
\frac{dP}{d\ga}\right]-E_P\left[\log\frac{dP_U}{d\mu}\right]\\
&=&E_\ga\left[\frac{dP}{d\ga} \log \frac{dP}{d\ga}\right]-E_{P_U}\left[\log\frac{dP_U}{d\mu}\right]\\
&=&E_\ga\left[\frac{dP}{d\ga} \log\frac{dP}{d\ga}\right]-E_\mu\left[\frac{dP_U}{d\mu} \log\frac{dP_U}{d\mu}\right]\\
&=&E_\ga[L\log L]-E\ga\left[\log \frac{dP_U}{d\mu}\circ U\right]\\
&=&\frac{1}{2}E_\ga\int_0^1|E_\mu[\dot{u}_s(w,m)|\calU_s(m)]|^2ds-E_\ga[-\log
\hat{\rho}]
\eeaa
and inserting the value of $\hat{\rho}$ given by the relation (\ref{ro})  completes the proof.

\qed

\remark
The similar results (cf. \cite{D,KZZ,M-Z}) in the literature concern the case where the
observation  $w\to U(w,m)$ is invertible $\ga$-almost surely,
consequently the first term is  reduced just to the half of the
$L^2(\mu,H)$-norm of $u$ (cf. \cite{ASU-3}).

\noindent
The following is a consequence of Bayes' lemma:

\begin{lemma}
\label{bayes1-lemma}
For any positive, measurable function $g$ on $W\times M$, we have 
$$
E_\ga[g|U]=\frac{1}{\hat{L}\circ
  U}\left(\int_ML(x,m)E_\mu\Big[g\mid U(\cdot,m)=x\Big]d\nu(m)\right)_{x=U}
$$
$\ga$-almost surely. In particular 
$$
E_\ga[g|U=x]=\frac{1}{\hat{L}(x)}\int_M L(x,m)E_\mu\Big[g\mid U(\cdot,m)=x\Big]d\nu(m)
$$
$P_U$ and $\mu$-almost surely.
\end{lemma}
\proof
Let $f\in C_b(W)$ and let $g$ be a positive, measurable function on
$W\times M$. We have 
\beaa
E_\ga[g\,f\circ U]&=&\int_M E_\mu[E_\mu[g\mid U(\cdot,m)]\,f\circ U(\cdot,m)]d\nu(m)\\
&=&\int_M \int_W L(w,m)\,E_\mu[g\mid U(\cdot,m)=w]\,f(w)d\mu(w)d\nu(m)\\
&=&\int_W f(w)\left(\int_M L(w,m) E_\mu[g\mid U(\cdot,m)=w]\,d\nu(m)\right)d\mu\\
&=&\int_W \frac{\hat{L}(w)}{\hat{L}(w)} f(w)\left(\int_M L(w,m) E_\mu[g\mid U(\cdot,m)=w]\,d\nu(m)\right)d\mu\\
&=&E_\ga\left[\frac{1}{\hat{L}\circ U}f\circ U\left(\int_M L(w,m)
    E_\mu[g\mid U(\cdot,m)=w]\,d\nu(m)\right)_{w=U}\right]
\eeaa
\qed

\noindent
From now on we return to the model $U_\la$ parametrized with $\la\in
\R$ and defined on the product space $W\times M$; namely
we assume that 
$$
U_\la(w,m)=w+u_\la(w,m)
$$
with the same independence hypothesis and the same regularity
hypothesis of $\la\to u_\la$ where the only difference consists of
replacement of the measure $\mu$ with the measure $\ga$ while
defining  the  spaces  $\DD_{p,k}$.

\begin{lemma}
\label{bayes2-lemma}
Let $\hat{L}_\la(w)$ denote the Radon-Nikodym derivative of $P_{U_\la}$
w.r. to $\mu$. We have 
$$
\hat{L}_\la(w)=\hat{L}_0(w)\exp \int_0^\la E_\ga\Big[\delta (K_\al
u'_\al)|U_\al=w\Big]d\al
$$
$\mu$-almost surely.
\end{lemma}
\proof
For any nice function $f$ on $W$, we have 
$$
\frac{d}{d\la}E_\ga[f\circ U_\la]=\frac{d}{d\la}E_\ga[f\,
L_\la]=\frac{d}{d\la}E_\mu[f\, \hat{L}_\la]\,.
$$
On the other hand 
\beaa
\frac{d}{d\la}E_\ga[f\circ U_\la]&=&E_\ga[f\circ U_\la \delta(K_\la
u_\la')]\\
&=&E_\ga[f\circ U_\la E_\ga[\delta(K_\la u_\la')|U_\la]]\\
&=&E_\ga[f L_\la(x,m) E_\ga[\delta(K_\la u_\la')|U_\la=x]]\\
&=&E_\mu[f \hat{L}_\la E_\ga [\delta(K_\la u_\la')|U_\la=x]]\,.
\eeaa

\qed

\remark
Note that we also have the following representation for $L_\la(w,m)$:
$$
L_\la(w,m)=L_0(w,m)\exp \int_0^\la E_\mu\Big[\delta (K_\al
u'_\al(\cdot,m))|U_\al(\cdot,m)=w\Big]d\al
$$
$\mu$-a.s.

\begin{lemma}
\label{deriv-lemma}
Let $\la\to \tau(\la)$ be defined as
$$
\tau(\la)=E_\ga[\hat{L}_\la\log \hat{L}_\la]\,,
$$
where
$\hat{L}_\la(w)=\int_ML_\la(w,m)d\nu(m)$ as before. We have
\beaa
\frac{d\tau(\la)}{d\la}&=&E_\ga\left[E_\ga[\delta(K_\la  u_\la')|U_\la]\log \hat{L}_\la\circ U_\la\right]\\
&=&E_\ga\left[E_\ga[\delta(K_\la  u_\la')|U_\la](-\log
  \hat{\rho}_\la)\right]
\eeaa
where $\hat{\rho}_\la$ is given by (\ref{ro}) as 
$$
\hat{\rho}_\la=\exp\left(-\int_0^1
  E_\ga[\dot{u}_\la(s)|\calU_\la(s)]dZ_\la(s)-\frac{1}{2}\int_0^1|E_\ga[\dot{u}_\la(s)|\calU_\la(s)]|^2ds\right)\,.
$$
Besides, we also have
$$
\frac{d^2\tau(\la)}{d\la^2}=E_\ga\left[E_\ga[\delta
  D_\la|U_\la](-\log \hat{\rho}_\la)+E_\ga[\delta(K_\la
  u'_\la)|U_\la]^2\right]
$$
where 
$$
D_\la=\delta(K_\la u'_\la)K_\la u'_\la+\frac{d}{d\la}K_\la u'_\la\,.
$$
\end{lemma}
\proof
The only thing that we need is the calculation of the second
derivative of $\hat{L}_\la$: let $f$ be a smooth function on $W$,
then, from Lemma \ref{bayes1-lemma},
\beaa
\frac{d^2}{d\la^2}E_\ga[f\circ U_\la]&=&\frac{d}{d\la}E_\ga[f\circ
U_\la\,\delta(K_\la u'_\la)]\\
&=&E_\ga\left[f\circ U_\ga\delta\left(\delta(K_\la u_\la')K_\la
  u_\la'+\frac{d}{d\la}(K_\la u_\la')\right)\right]\\
&=&E_\ga[f\circ U_\ga\,\delta D_\la]\\
&=&E_\ga[f(x)\,E_\ga[\delta D_\la|U_\la=x]\,\hat{L}_\la(x)]\,.
\eeaa
\qed

\noindent
As an immediate consequence we get
\begin{corollary}
We have the following relation:
\beaa
\frac{d^2}{d\la^2}I(U_\la,m)&=&E_\ga\Big[E_\mu[\delta
(D_\la(\cdot,m))|U_\la(m)]
(-\log E_\mu[\rho(-\delta
u_\la(\cdot,m))|U_\la(m)])\\
&&+E_\mu[\delta(K_\la u'_\la(\cdot,m))|U_\la(m)]^2\Big]\\
&&-E_\ga\Big[E_\ga[\delta
(D_\la)|U_\la](-\log \hat{\rho}_\la)+E_\ga[\delta(K_\la u'_\la)
|U_\la]^2\Big]\,.
\eeaa
\end{corollary}

\vspace{2cm}
\footnotesize{
\noindent
A. S. \"Ust\"unel, Institut Telecom, Telecom ParisTech, LTCI CNRS D\'ept. Infres, \\
46, rue Barrault, 75013, Paris, France\\
ustunel@enst.fr}

\end{document}